\documentclass[preprint,12pt]{elsarticle}
\usepackage{amsmath,url}

\newtheorem{theorem}{\textbf{Theorem}}

\newtheorem{remark}{\textbf{Remark}}

\newtheorem{lemma}{\textbf{Lemma}}

\usepackage{mathrsfs}

\def\a {\alpha}
\def\p {\alpha}

\def\Z {\mathbb{Z}}
\def\Q {\mathbb{Q}}

\usepackage{amsfonts}
\usepackage{amssymb}

\def\e {\epsilon}

\def\d {\delta}
\def\l {\lambda}
\def\dfrac {\displaystyle\frac}
\def\g {\gamma}

\def\tareesidedbox#1{\setbox0=\hbox{$#1$}\dimen0=\wd0 \advance\dimen0 by3pt\rlap{\hbox{\vrule height8pt width.4pt
 depth2pt \kern-.4pt\vrule height8.4pt width\dimen0 depth-8pt\kern-.4pt \vrule height8pt width.4pt depth2pt}}
 \relax \hbox to\dimen0{\hss$#1$\hss}}

\usepackage{amssymb}

\journal{...}

\begin{document}

\begin{frontmatter}



\title{Generalized Cullen Numbers in Linear Recurrence Sequences}



\author{Yuri Bilu\fnref{ybb}}{\ead{yuri@math.u-bordeaux.fr}}
\author[dmm]{Diego Marques\fnref{dm}\corref{CA}}{\ead{diego@mat.unb.br}}
\author[at]{Alain Togb\'e \fnref{att}}{\ead{atogbe@pnw.edu}}
\address[ybb]{IMB, Universit\' e de Bordeaux 1,
351 cours de la Lib\' eration, 33405, Talence CEDEX, France}
\address[dmm]{Departamento de Matem\' atica, Universidade  de Bras\' ilia, Bras\' ilia, 70910-900, Brazil}
\address[at]{Department of Mathematics, Statistics and Computer Science, Purdue University Northwest, 1401 S, U.S. 421, Westville, IN 46391, USA}


\fntext[dm]{Supported by DPP/UnB and a postdoctoral fellowship (247311/2012-0) from CNPq-Brazil}
\fntext[att]{Supported in part by Purdue University North Central}

\cortext[CA]{Corresponding author}

\begin{abstract}
A Cullen number is a number of the form $m2^m+1$, where $m$ is a positive integer. In 2004, Luca and St\u anic\u a proved, among other things, that the largest Fibonacci number in the Cullen sequence is $F_4=3$. Actually, they searched for generalized Cullen numbers among some binary recurrence sequences. In this paper, we will work on higher order recurrence sequences. For a given linear recurrence $(G_n)_n$, under weak assumptions, and a given polynomial $T(x)\in \Z[x]$, we shall prove that if $G_n=mx^m+T(x)$, then
\[
m\ll\log \log |x|\log^2(\log \log |x|)\ \mbox{and}\ n\ll\log |x|\log\log |x|\log^2(\log \log |x|),
\]
where the implied constant depends only on $(G_n)_n$ and $T(x)$.
\end{abstract}


\end{frontmatter}



\section{Introduction}

A {\it Cullen number} is a number of the form $m2^m+1$ (denoted by $C_m$), where $m$ is a nonnegative integer. A few terms of this sequence are
\[
1, 3, 9, 25, 65, 161, 385, 897, 2049, 4609, 10241, 22529,\ldots
\]
which is the OEIS \cite{oeis} sequence A002064 (this sequence was introduced in 1905 by the Father J. Cullen \cite{cul} and it was mentioned in the well-known Guy's book \cite[Section {\bf B20}]{guy}). These numbers gained great interest in 1976, when C. Hooley \cite{hoo} showed that almost all Cullen numbers are composite. However, despite being very scarce, it is still conjectured the existence of infinitely many \textit{Cullen primes}. For instance,  $C_{6679881}$ is a prime number with more than 2 millions of digits (PrimeGrid, August 2009).

These numbers can be generalized to the \textit{generalized Cullen numbers} which are numbers of the form
\begin{center}
$C_{m,s}=ms^m+1$,
\end{center}
where $m\geq 1$ and $s\geq 2$. Clearly, one has that $C_{m,2}=C_m$, for all $m\geq 1$. For simplicity, we call $C_{m,s}$ of \textit{$s$-Cullen number}. This family was introduced by H. Dubner \cite{dub} and is one of the main sources for prime number ``hunters". A big prime of the form $C_{m,s}$ is $C_{139948,151}$ an integer with 304949 digits.

Many authors have searched for special properties of Cullen numbers and their generalizations. Concerning these numbers, we refer to \cite{pt,uber,new} for primality results and \cite{mdc} for their greatest common divisor. The problem of finding Cullen numbers belonging to others known sequences has attracted much attention in the last two decades. We cite \cite{pseudo} for pseudoprime Cullen numbers, and \cite{bar} for Cullen numbers which are both Riesel and Sierpi\' nski numbers.

A sequence $(G_n)_{n\geq 0}$ is \textit{a linear recurrence sequence} with coefficients $c_0$, $c_1,$\ldots,$c_{k-1}$, with $c_0\neq 0,$ if
\begin{equation}\label{recu}
G_{n+k}=c_{k-1}G_{n+k-1}+\cdots + c_1G_{n+1}+c_0G_n,
\end{equation}
for all positive integer $n$. A recurrence sequence is therefore completely determined by the \textit{initial values} $G_0,\ldots,G_{k-1}$, and by the coefficients $c_0,c_1,\ldots,c_{k-1}$. The integer $k$ is called the {\it order} of the linear recurrence. The \textit{characteristic polynomial} of the sequence $(G_n)_{n\geq 0}$ is given by
$$G(x)=x^{k}-c_{k-1}x^{k-1}-\cdots - c_1x-c_0.$$
It is well-known that for all $n$
\begin{equation}\label{general}
G_n=g_1(n)r_1^n+\cdots +g_{\ell}(n)r_{\ell}^n,
\end{equation}
where $r_j$ is a root of $G(x)$ and $g_j(x)$ is a polynomial over a certain number field, for $j=1,\ldots,\ell$. In this paper, we consider only integer recurrence sequences, i.e., recurrence sequences whose coefficients and initial values are integers. Hence, $g_j(n)$ is an algebraic number, for all $j=1,\ldots,{\ell}$, and $n\in \mathbb{Z}$.

A general Lucas sequence $(U_n)_{n\geq 0}$ given by $U_{n+2}=aU_{n+1}+bU_n$, for $n\geq 0$, where the values $a,\ b,\ U_0$ and $U_1$ are previously fixed, is an example of a linear recurrence of order $2$ (also called {\it binary}). For instance, if $U_0=0$ and $U_1=a=b=1$, then $(U_n)_{n\geq 0}=(F_n)_{n\geq 0}$ is the well-known \textit{Fibonacci sequence} and for $U_0=2$ and $U_1=a=b=1$, then $(U_n)_{n\geq 0}=(L_n)_{n\geq 0}$ is the sequence of the {\it Lucas numbers}:

In 2003, Luca and St\u anic\u a \cite{LS} showed, in particular, that Cullen numbers occur only finitely many times in a binary recurrent sequence satisfying some additional conditions. As application, they proved that the largest Fibonacci number in the Cullen sequence is $F_4=3=1\cdot 2^1+1$. Very recently, Marques \cite{DM} searched for Fibonacci numbers in $s$-Cullen sequences. In particular, he proved that there is no Fibonacci number that is also a nontrivial $s$-Cullen number when all divisors of $s$ are not Wall-Sun-Sun primes (i.e., $p^2\nmid F_{p-(5/p)}$). See also \cite{APDM}. We remark that no Wall-Sun-Sun prime is known as of July 2017, moreover if any exist, they must be greater than $2.6\cdot 10^{17}$.

In this paper, we are interested in much more general Cullen numbers among terms of linear recurrences. More precisely, our goal is to work on the Diophantine equation
\begin{equation}\label{M}
G_n=mx^m+T(x),
\end{equation}
for a given polynomial $T(x)\in \Z[x]$. Observe that when $x$ is fixed and $T(x)= 1$ the right-hand side of (\ref{M}) is an $x$-Cullen number.

We remark that several authors investigated the related equation
\[
G_n=ax^m+T(x),
\]
where $a$ is fixed. Among the results for a general recurrence (under some technical hypotheses), it was proved finiteness of solutions for $T(x)=0$ by Shorey and Stewart \cite{SS}, for $T(x)=c$ by Stewart \cite{Ste} and for any $T(x)$ by Nemes and Peth\H o \cite{Nemes}. Moreover, these results are effective.

Here, our main result is the following

\begin{theorem}\label{main1}
Let $(G_n)_n$ be an integer linear recurrence with roots $r_1,\ldots,r_k$ satisfying either
\begin{itemize}
\item[(i)] $|r_1|>1>|r_j|>0$, for $j>1$; or
\item[(ii)] $|r_1|>|r_2|>|r_j|>0$, for $j>2$.
\end{itemize}
Moreover, we suppose that $r_1$ is a simple root. Let $T(x)\in \Z[x]$ be a polynomial. There exist effectively computable constants $C_1, C_2$, depending only on $(G_n)_n$ and $T(x)$, such that if $(m,n,x)$ is a solution of the Diophantine equation (\ref{M}), then
\[
m\leq C_1\log \log |x|\log^2(\log \log |x|)\ \mbox{and}\ n\leq C_2\log |x|\log\log |x|\log^2(\log \log |x|).
\]
\end{theorem}

Observe that we cannot ensure here finitely many values for $|x|$. For example, for $G_n=2L_n$ and $T(x)=-4$, one has that $(n,m,x)=(4t,2,L_{2t})$ is solution for Eq. (\ref{M}) for all $t\geq 0$.

Let $k\geq 2$ and denote $F^{(k)}:=(F_n^{(k)})_{n\geq -(k-2)}$, the \textit{$k$-generalized Fibonacci sequence} whose terms satisfy the recurrence relation
\begin{equation}\label{rec}
F_{n+k}^{(k)}=F_{n+k-1}^{(k)}+F_{n+k-2}^{(k)}+\cdots + F_{n}^{(k)},
\end{equation}
with initial conditions $0,0,\ldots,0,1$ ($k$ terms) and such that the first nonzero term is $F_1^{(k)}=1$.

The above sequence is one among the several generalizations of Fibonacci numbers. Such a sequence is also called $k$-\textit{step Fibonacci sequence}, the \textit{Fibonacci $k$-sequence}, or $k$-\textit{bonacci sequence}. Clearly for $k=2$, we obtain the classical Fibonacci numbers, for $k=3$, the Tribonacci numbers, for $k=4$, the Tetranacci numbers, etc.

Recently, these sequences have been the main subject of many papers. We refer to \cite{BLP} for results on the largest prime factor of $F_n^{(k)}$ and we refer to \cite{BL2} for the solution of the problem of finding powers of two belonging to these sequences. In 2013, two conjectures concerning these numbers were proved. The first one, proved by Bravo and Luca \cite{BL} is related to \textit{repdigits} (i.e., numbers with only one distinct digit in its decimal expansion) among $k$-Fibonacci numbers (proposed by Marques \cite{util}) and the second one, a conjecture (proposed by Noe and Post \cite{noe}) about coincidences between terms of these sequences, proved independently by Bravo-Luca \cite{BL0} and Marques \cite{Marques0}.

If we use Theorem \ref{main1} to sequence $(G_n)_n=(F_n^{(k)})_n$, we get finitely many solutions for Eq. (\ref{M}), for each  $k\geq 2$. However, we shall improve the method and we find an upper bound for the number of Cullen numbers (case $x=2$ and $T(x)=1$) in $\cup_{k\geq 2}F^{(k)}$. More precisely,
\begin{theorem}\label{main2}
If $(m,n,k)$ is a solution of the Diophantine equation
\begin{equation}\label{Main2}
F_n^{(k)}=m2^m+1
\end{equation}
in positive integers $m,n$ and $k\geq 2$, then
\begin{center}
$m<9.5\cdot 10^{23}, n<2.4\cdot 10^{24}$ and $k\leq 158$.
\end{center}
\end{theorem}

Let us give a brief overview of our strategy for proving Theorem \ref{main2}. First, we use a Dresden and Du formula \cite[Formula (2)]{dres} to get an upper bound for a linear form in three logarithms related to equation (\ref{Main2}). After, we use a lower bound due to Matveev to obtain an upper bound for $m$ and $n$ in terms of $k$. Very recently, Bravo and Luca solved the equation $F_n^{(k)}=2^m$ and for that they used a nice argument combining some estimates together with the Mean Value Theorem (this can be seen in pages $77$ and $78$ of \cite{BL2}). In our case, we use Bravo-Luca's approach to get an inequality involving a linear form in two logarithms. In the other case, we use a lower bound due to Laurent to get substantially upper bounds for $m, n$ and $k$.  The computations in the paper were performed using \textit{Mathematica}${}^\text{\tiny\textregistered}$

\section{Auxiliary results}\label{sec2}

In this section, we recall some results that will be very useful for the proof of the above theorems. Let $G(x)$ be the characteristic polynomial of a linear recurrence $G_n$. One can factor $G(x)$ over the set of complex numbers as
$$
G(x)=(x-r_1)^{m_1}(x-r_2)^{m_2}\cdots (x-r_{\ell})^{m_{\ell}},
$$
where $r_1,\ldots,r_{\ell}$ are distinct non-zero complex numbers (called the {\it roots} of the recurrence) and $m_1,\ldots,m_{\ell}$ are positive integers.  A root $r_j$ of the recurrence is called a {\it dominant root} if $|r_j|>|r_i|$, for all $j\neq i\in \{1,\ldots,\ell\}$. The corresponding polynomial $g_j(n)$ is named the \textit{dominant polynomial} of the recurrence. A fundamental result in the theory of recurrence sequences asserts that there exist uniquely determined non-zero polynomials $g_1,\ldots,g_{\ell}\in \mathbb{Q}(\{r_j\}_{j=1}^{\ell})[x]$, with $\deg g_j\leq m_j-1$, for $j=1,\ldots,\ell$, such that
\begin{equation}\label{written}
G_n=g_1(n)r_1^n+\cdots +g_{\ell}(n)r_{\ell}^n,\ \mbox{for\ all}\ n.
\end{equation}
For more details, see \cite[Theorem C.1]{shorey}.

In the case of the Fibonacci sequence, the above formula is known as {\it Binet's formula}:
\begin{equation}\label{binet1}
F_n=\displaystyle\frac{\alpha^n-\beta^n}{\alpha-\beta},
\end{equation}
where $\alpha=(1+\sqrt{5})/2$ (the golden number) and $\beta=(1-\sqrt{5})/2=-1/\alpha$.
Equation (\ref{written}) and some tricks will allow us to obtain linear forms in three logarithms and then determine lower bounds {\it \`a la Baker} for these linear forms. From the main result of Matveev \cite{matveev}, we deduce the following lemma.

\begin{lemma}\label{lemma1}
Let $\gamma_1,\ldots,\gamma_t$ be real algebraic numbers and let $b_1,\ldots,b_t$ be nonzero rational integer numbers. Let $D$ be the degree of the number field $\Q(\gamma_1,\ldots,\gamma_t)$ over $\Q$ and let $A_j$ be a positive real number satisfying
\begin{center}
$A_j\geq \max\{Dh(\g_j),|\log \g_j|,0.16\}$, for $j=1,\ldots,t$.
\end{center}
Assume that
$$
B\geq \max\{|b_1|,\ldots,|b_t|\}.
$$
If $\g_1^{b_1}\cdots \g_t^{b_t}\neq 1$, then
$$
|\g_1^{b_1}\cdots \g_t^{b_t}-1|\geq \exp(-1.4\cdot 30^{t+3}\cdot t^{4.5}\cdot D^2(1+\log D)(1+\log B)A_1\cdots A_t).
$$
\end{lemma}
As usual, in the previous statement, the \textit{logarithmic height} of an $n$-degree algebraic number $\alpha$ is defined as
$$
h(\alpha)=\frac{1}{n}(\log |a|+\sum_{j=1}^n\log \max\{1,|\alpha^{(j)}|\}),
$$
where $a$ is the leading coefficient of the minimal polynomial of $\alpha$ (over $\mathbb{Z}$) and $(\alpha^{(j)})_{1\leq j\leq n}$ are the conjugates of $\alpha$.


Now, we are ready to deal with the proofs of our results. 

\section{The proof of Theorem \ref{main1}} \label{sec3}

Throughout the proof, the numerical constants implied by $\ll$ depend only on $(G_n)_n$ and $T(x)$. Also, without loss of generality, we may suppose $|x|\geq 2$ (i.e., $x\neq -1,0,1$).

First, since $r_1$ is a simple dominant root then $g_1(n)$ in formula (\ref{written}) is a constant, say $g$ (because the degree of $g_1(n)$ would be at most $m_1-1=1-1=0$). Now, we rewrite Eq. (\ref{M}) as
\[
mx^m-gr_1^n=B(x,n),
\]
where $B(x,n):=\sum_{j=2}^kg_j(n)r_j^n-T(x)$. Then
\[
\dfrac{mx^mr_1^{-n}}{g}-1=\dfrac{B(x,n)}{gr_1^n}.
\]
If $B(x,n)=0$ and (ii) holds, we use the same argument than Nemes and Peth\H o to get $m\ll 1$ and the proof is complete (see lines 25-35 in page 231 of \cite{Nemes}). However, if $B(x,n)=0$ and (i) holds, we get the relation $mx^m=gr_1^n$. So, we can take the conjugates of this relation in $\Q(r_1)$ to get $mx^m=g^{(t)}r_t^n$, where the $g_1^{(i)}$'s are the conjugates of $g_1$ over $\Q(r_1)$. Thus, by taking absolute values and using that $|r_t|<1$ we obtain $m2^m\ll 1$ yielding $m\ll 1$. 

Thus, we may suppose $B(x,n)\neq 0$ and in this case Nemes and Peth\H o  \cite[p. 232]{Nemes} proved that $|B(x,n)|\leq r_1^{n(1-\delta)}$, for some $\delta\ll 1$. Therefore

\begin{equation}
\left|\dfrac{mx^mr_1^{-n}}{g}-1\right|\ll \dfrac{1}{r_1^{n\d}}.
\end{equation}

Let $\Lambda=\log (m/g)-n\log r_1+m\log x$. Since $x<e^x-1$ and for $x<0$, $|e^x-1|=1-e^{-|x|}$, then the previous inequality yields $|\Lambda|\ll 1/r_1^{n\delta+O(1)}$ yielding

\begin{equation}\label{est1}
\log |\Lambda| \ll -(n\delta+O(1))\log r_1.
\end{equation}

Now, we will apply Lemma \ref{lemma1}. To this end, take
\[
t:=3,\;\; \g_1:=m/g,\ \g_2:=x,\ \g_3:=r_1,
\]
and
\[
b_1:=1,\ b_2:= m,\ b_3:= -n.
\]
For this choice, we have $D=[\Q(g,r_1):\Q]\leq k!$. Also $h(\g_1)\leq \log m+h(g)\ll \log m$, $h(\g_2)=\log |x|$ and $h(\g_3)\ll \log r_1$, where we used the well-known facts that $h(xy)\leq h(x)+h(y)$ and $h(x)=h(x^{-1})$. 

Note that Eq. (\ref{M}) implies that $m\log |x|\asymp n$. In fact, one has that
\[
r_1^{n-O(1)}\ll |G_n|=|mx^m+T(x)|\ll mx^m,
\]
where we used that $|T(x)|\ll |x|^{\deg T}$. Thus, we obtain $n\ll m\log |x|$. On the other hand, $m|x|^m\leq |g|r_1^n+|B(x,n)|\ll r_1^{O(n)}$ (here we used that $|B(x,n)|\leq r_1^{n(1-\d)}$). By applying the $\log$ function we arrive at
\[
1+m\log |x|\ll \log m+m\log |x|\ll n
\]
and thus $m\log |x|\ll n$. Therefore, we have that $B \ll m\log |x|$.

Since $B(x,n)\neq 0$, the left-hand side of (\ref{est1}) is nonzero and so the conditions to apply Lemma \ref{lemma1} are fulfilled yielding
\begin{equation}\label{est2}
\log |\Lambda|>-\log^2 m\log |x|\log \log |x|.
\end{equation}
Combining estimates (\ref{est1}) and (\ref{est2}) we have
\begin{equation}\label{A}
n\ll \log^2 m\log |x|\log \log |x|.
\end{equation}

Combining this estimate with (\ref{A}) we get
\[
\dfrac{m}{\log ^2 m}\ll \log \log |x|.
\]
As we shall prove in (\ref{key}), the inequality above implies
\[
m\ll \log \log |x|\log^2 (\log \log |x|).
\]

Now, we use the estimate $n\ll m\log |x|$ to get the desired inequality on $n$, i.e.,
\[
n\ll \log |x|\log \log |x|\log^2 (\log \log |x|).
\]
The proof is then complete.
\qed

\section{The proof of Theorem \ref{main2}}

\subsection{Auxiliary results}

Before proceeding further, we will recall some facts and properties of these sequences which will be used
after.

We know that the characteristic polynomial of $(F_n^{(k)})_n$ is
\[
\psi_k(x):=x^k-x^{k-1}-\cdots -x-1
\]
and it is irreducible over $\Q[x]$ with just one zero outside the unit circle. That single zero is located between $2(1-2^{-k})$ and $2$ (as it can be seen in \cite{HUA}). Also, in a recent paper, G. Dresden and Z. Du \cite[Theorem 1]{dres} gave a simplified ``Binet-like"\ formula for $F_n^{(k)}$:
\begin{equation}\label{binet}
F_n^{(k)}=\displaystyle\sum_{i=1}^k\dfrac{\alpha_i-1}{2+(k+1)(\a_i-2)}\a_i^{n-1},
\end{equation}
for $\alpha=\alpha_1,\ldots,\a_k$ being the roots of $\psi_k(x)$. Also, it was proved in \cite[Lemma 1]{BL} that
\begin{equation}\label{Lu}
\alpha^{n-2}\leq F_n^{(k)}\leq \alpha^{n-1},\ \mbox{for\ all}\ n\geq 1,
\end{equation}
where $\a$ is the dominant root of $\psi_k(x)$. Also, the contribution of the roots inside the unit circle in formula (\ref{binet}) is almost trivial. More precisely, it was proved in \cite{dres} that
\begin{equation}\label{small}
|F_n^{(k)}-g(\alpha,k)\alpha^{n-1}|<\dfrac{1}{2},
\end{equation}
where we adopt throughout the notation $g(x,y):=(x-1)/(2+(y+1)(x-2))$.

Very recently, Bravo and Luca \cite{BL2} found all powers of two in $k$-generalized Fibonacci sequences. Their nice method can be slightly changed to show that $(n,k,m)=(1,4,2)$ and $(5,2,2)$ are the only solutions of the equation $F_n^{(k)}=2^m+1$, with $k\geq 2$. Thus, the only solution of Eq. (\ref{Main2}) such that $m$ is a power of two is $(n,k,m)=(1,4,2)$. So, throughout the paper, we shall suppose that $m$ is not a power of two and that $m\geq 10$ (the case $m<10$ can be easily solved). Note also that, by definition, $F_n^{(k)}$ is a power of two for all $1\leq n\leq k+1$ and hence these values cannot be Cullen numbers. Thus, it is enough to consider $n>k+1$. Finally, due to \cite[Theorem 3]{LS}, we can suppose that $k\geq 3$.

\subsection{The proof}

First, we use Eq. (\ref{Main2}) together with the formula (\ref{binet}) to obtain
\begin{equation}\label{s}
g(\a,k)\a^{n-1}-m2^m=1-\displaystyle\sum_{i=2}^kg(\a_i,k)\a_i^{n-1}\in (1/2,3/2),
\end{equation}
 where we used (\ref{small}). Thus, equation (\ref{s}) implies that
$$0<g(\a,k)\a^{n-1}-m2^m <3/2.$$
So, dividing by $m2^m$, we get
\begin{equation}\label{Est1}
\left|\dfrac{g(\a,k)\a^{n-1}}{m2^m}-1\right| < 1/2^{m+1},
\end{equation}
for $ m\geq 3.$

In order to use Lemma \ref{lemma1}, we take
\[
t:=3,\;\; \g_1:=g(\a,k)/m,\ \g_2:=2,\ \g_3:=\a
\]
and
\[
b_1:=1,\ b_2:= -m,\ b_3:= n-1.
\]
For this choice, we have $D=[\Q(\a):\Q]= k$. Also $h(\g_1)\leq \log ((4k+4)m)$, $h(\g_2)=\log 2$ and $h(\g_3)< 0.7/k$. Thus, we can take $A_1:= k\log ((4k+4)m), A_2:= k\log 2$ and $A_3:=0.7$.

Moreover, using the inequalities (\ref{Lu}), we get
\[ (7/4)^{n-2}<\a^{n-2}<F_n^{(k)}=m2^m+1<2^{2m-1}\]
and so $n<2.5m+0.8$. Note that $\max\{|b_1|,|b_2|,|b_3|\}=\max\{m,n-1\}\leq 2.5m+0.8=:B$. Since $g(\a,k)\a^{n-1}2^{-m}/m> 1$ (by (\ref{s})), we are in position to apply Lemma \ref{lemma1}. This lemma together with a straightforward calculation gives
\begin{equation}\label{Est2}
\left|\dfrac{g(\a,k)\a^{n-1}}{m2^m}-1\right| >  \exp(-6.7\cdot 10^{11}k^{4}\log^2 m),
\end{equation}
where we used that $1+\log k<2\log k$, for $k\geq 2$, $1+\log (2.5m+0.8)<1.9\log m$, for $m\geq 10$, and $\log ((4k+4)m)<2.5\log m$ (to prove this last inequality, we used that $2.5m+0.8>n>k+1$).

By combining (\ref{Est1}) and (\ref{Est2}), we obtain
\[
\dfrac{m}{\log^2 m}<9.7\cdot 10^{11}k^{4}\log k.
\]
Since the function $x\mapsto x/\log^2 x$ is increasing for $x>e$, then it is a simple matter to prove that
\begin{equation}\label{key}
\dfrac{x}{\log^2 x}<A\ \ \mbox{implies\ that}\ \ x<2A\log^2 A\ (\mbox{for}\ A\geq 10^7).
\end{equation}
In fact, suppose the contrary, i.e. $x\geq 2A\log^2 A$. Then
\[
\dfrac{x}{\log^2 x}\geq \dfrac{2A\log^2 A}{\log^2 (2A\log^2 A)}>A,
\]
which contradicts our inequality. Here we used that $\log^2 (2A\log^2 A)<2\log^2 A$, for $A\geq 10^7$.

Thus, using (\ref{key}) for $x:=m$ and $A:=9.7\cdot 10^{11}k^{4}\log k$, we have that
\[
m<2(9.7\cdot 10^{11}k^{4}\log k)\log^2 (9.7\cdot 10^{11}k^{4}\log k).
\]
A straightforward calculation gives
\begin{equation}\label{m<k}
m<5.9\cdot 10^{13}k^{4}\log^2 k.
\end{equation}

Now, we shall prove that there is no solution when $k\geq 159$. In this case, (\ref{m<k}) implies
\[
n<2.5m+0.8<11.8\cdot 10^{13}k^{4}\log^2 k + 0.8<2^{k/2}.
\]

Now, we use a key argument due to Bravo and Luca \cite[p. 77-78]{BL2}. 

Setting $\l=2-\p$, we deduce that $0<\l<1/2^{k-1}$ (because $2(1-2^{-k})<\p<2$). So
\[
\p^{n-1}=(2-\l)^{n-1}=2^{n-1}\left(1-\dfrac{\l}{2}\right)^{n-1}>2^{n-1}(1-(n-1)\l),
\]
since that the inequality $(1-x)^n>1-2nx$ holds for all $n\geq 1$ and $0<x<1$. Moreover, $(n-1)\l<2^{k/2}/2^{k-1}=2/2^{k/2}$ and hence
\[
2^{n-1}-\dfrac{2^n}{2^{k/2}}<\p^{n-1}<2^{n-1}+\dfrac{2^n}{2^{k/2}},
\]
yielding
\begin{equation}\label{delta}
|\p^{n-1}-2^{n-1}|<\dfrac{2^n}{2^{k/2}}.
\end{equation}

Now, we define for $x>2(1-2^{-k})$ the function $f(x):=g(x,k)$ which is differentiable in the interval $[\p,2]$. So, by the Mean Value Theorem, there exists $\xi\in (\p,2)$, such that $f(\p)-f(2) = f'(\xi)(\p-2)$. Thus
\begin{equation}\label{eta}
|f(\p)-f(2)|  < \dfrac{2k}{2^{k}}, 
\end{equation}
where we used the bounds $|\p-2|<1/2^{k-1}$ and $|f'(\xi)|<k$. For simplicity, we denote $\d=\p^{n-1}-2^{n-1}$ and $\eta=f(\p)-f(2)=f(\p)-1/2$. After some calculations, we arrive at
\[
2^{n-2}=f(\p)\p^{n-1}-2^{n-1}\eta - \dfrac{\d}{2} - \d \eta.
\]
Therefore
\begin{eqnarray*}
|2^{n-2}-m2^m| & \leq & \dfrac{3}{2}+2^{n-1}|\eta| + \left|\dfrac{\d}{2}\right| + |\d \eta |\\
 & \leq & \dfrac{3}{2}+\dfrac{2^nk}{2^{k}}+\dfrac{2^{n-1}}{2^{k/2}} + \dfrac{2^{n+1}k}{2^{3k/2}},
\end{eqnarray*}
where we used (\ref{delta}) and (\ref{eta}). Since $n>k+1$, one has that $2^{n-2}/2^{k/2}\geq 2^{k/2}>3/2$ (for $k\geq 2$) and we rewrite the above inequality as
\[
|2^{n-2}-m2^m|<\dfrac{2^{n-2}}{2^{k/2}} + \left(\dfrac{4k}{2^{k/2}}\right)\dfrac{2^{n-2}}{2^{k/2}}+ 2\cdot \dfrac{2^{n-2}}{2^{k/2}}+ \left(\dfrac{8k}{2^{k}}\right)\dfrac{2^{n-2}}{2^{k/2}}.
\]
Since the inequality $\max_{k\geq 159}\{4k/2^{k-1}, 8k/2^{k/2}\}<5.4\cdot 10^{-22}$ holds, then
\begin{equation}\label{Phi}
|2^{n-2}-m2^{m}|<\dfrac{3.2\cdot 2^{n-2}}{2^{k/2}},
\end{equation}
or equivalently
\begin{equation}\label{ss}
|1-m2^{-(n-m-2)}|<\dfrac{3.2}{2^{k/2}}.
\end{equation}

Since $m\geq 10$, we have
\begin{itemize}
\item If $\log m/\log 2+m+3\leq n$, then $1-m/2^{n-m-2}\geq 1/2$ yielding $k\leq 5$;
\item If $\log m/\log 2+m+1\geq n$, then $m2^{n-m-2}-1\geq 1$ leading to $k\leq 3$
\end{itemize}
which is not possible. Since $\log m/\log 2\notin \Q$ when $m$ is not a power of $2$, we may suppose that $n= \lfloor \log m/\log 2\rfloor +m+\d$, for $\d\in\{2,3\}$.

Note that (\ref{ss}) is equivalent to
\begin{equation}
|1-e^{\Lambda}|<\dfrac{3.2}{2^{k/5}},
\end{equation}
where $\Lambda:=\log m-(\lfloor \log m/\log 2\rfloor+\d-2)\log 2$.

Since $m$ is not a power of $2$, then $m$ and $2$ are multiplicatively independent. In particular, $\Lambda\neq 0$. If $\Lambda>0$, then $\Lambda<e^{\Lambda}-1<3.2/2^{k/2}$. In the case of $\Lambda<0$, we use $1-e^{-|\Lambda|}=|e^{\Lambda}-1|<3.2/2^{k/2}$ to get $e^{|\Lambda|}<1/(1-3.2\cdot 2^{-k/2})$. Thus
\[
|\Lambda|<e^{|\Lambda|}-1<\dfrac{3.2\cdot 2^{-k/2}}{1-3.2\cdot 2^{-k/2}}<3.6\cdot 2^{-k/2},
\]
where we used that $1/(1-3.2\cdot 2^{-k/2})<1.1$, for $k\geq 159$. In any case, we have
\begin{equation}\label{sq}
|\Lambda|<3.6\cdot 2^{-k/2}
\end{equation}
and so
\begin{equation}\label{est3}
\log |\Lambda|<\log (3.6)-\frac{k}{2}\log 2.
\end{equation}

Now, we will determine a lower bound for $\Lambda$. We remark that the bounds available for linear forms in two logarithms are substantially better than those available for linear forms in three logarithms. Here we choose to use a result due to Laurent \cite[Corollary 2]{Laurent} with $m=24$ and $C_2=18.8$. First let us introduce some notations. Let $\alpha_1,\a_2$ be real algebraic numbers, with $|\a_j|\geq 1$, $b_1,b_2$ be positive integer numbers and $$\Lambda=b_2\log \a_2-b_1\log \a_1.$$

Let $A_j$ be real numbers such that
$$
\log A_j\geq \max\{h(\a_j),|\log \a_j|/D, 1/D\}, j\in \{1,2\},
$$
where $D$ is the degree of the number field $\Q(\alpha_1,\alpha_2)$ over $\Q$. Define
$$
b'=\dfrac{b_1}{D\log A_2}+\dfrac{b_2}{D\log A_1}.
$$
Laurent's result asserts that if $\a_1,\a_2$ are multiplicatively independent, then
\[
\log |\Lambda| \geq  -18.8\cdot D^4\left(\max\{\log b'+0.38,m/D,1\}\right)^2 \cdot \log A_1\log A_2.
\]

We then take
\[
D=1,\; b_1=\lfloor \log m/\log 2\rfloor+\d-2,\; b_2=1,\; \a_1=2,\; \a_2=m.
\]
We choose $\log A_1=1$ and $\log A_2=\log m$. So we get
$$
b'=\dfrac{\lfloor \log m/\log 2\rfloor+\d-2}{\log m}+1< \dfrac{1}{\log m}+\dfrac{1}{\log 2}+1.
$$

Thus, by Corollary $2$ of \cite{Laurent} we get
\begin{equation}\label{est4}
\log |\Lambda|\geq -13.1\cdot 24^2\log m.
\end{equation}
Now, we combine the estimates (\ref{est3}) and (\ref{est4}) to obtain $k<21761\log m$. On the other hand, inequality (\ref{m<k}) gives $k<2.1\cdot 10^6$. Therefore, $m<2.5\cdot 10^{41}$. Now, we come back to (\ref{sq}) and by using Mathematica, we arrive at
\begin{eqnarray*}
\frac{2^{-k/2}}{\log 2} & > & \displaystyle\min_{\theta\in\{0,1\},3\leq m\leq 2.5 \cdot 10^{41}, m\neq 2^s}\left\{\left |\frac{\log m}{\log 2}-\left(\left \lfloor \frac{\log m}{\log 2}\right \rfloor+\theta\right)\right |\right\}\\
& \geq & \min_{3\leq m\leq 2.5 \cdot 10^{41}, m\neq 2^s}\min\left\{\left\{\frac{\log m}{\log 2}\right\}, 1- \left\{\frac{\log m}{\log 2}\right\}\right \}\\
& > & 8.2\cdot 10^{-42},
\end{eqnarray*}
where this minimum occurs when $m=2^{137}\pm 1$ (here, as usual $\{x\}$ denotes the fractional part of a real number $x$). This yields $k\leq 274$ implying $m<1.1\cdot 10^{25}$. Now, we repeat the above process two times (with the minimum occuring in $m=2^{83}\pm 1$ and $m=2^{79}\pm 1$) to obtain $k\leq 158$. This contradicts the assumption of $k\geq 159$.\qed

\begin{remark}
We remark to the reader that it must be possible to improve the upper bound for $m$, $n$ and $k$ in Theorem \ref{main2}. Unfortunately, it is not possible to decrease them to fulfill all remaining cases. On the other hand, the usual approach to finish the finite many cases is by using the Baker-Davenport reduction method (mainly, results related to a Dujella-Peth\"{o} theorem). However, for this problem, we have a form like
\[
(n-1)\gamma_k-m+\mu_{k,m},
\]
where $\gamma_k:=\log \a/\log 2$ and $\mu_{m,k}:=\log (g(\a,k)/m)/\log 2$. To use the reduction method, we should get a positive lower bound for a quantity (called $\e$) depending on, in this case, $k$ and $m$. The problem here is the dependence on $m$ which by its size ($\approx 10^{23}$) becomes the calculation ``impossible", by our computational tools.
\end{remark}

\section*{Acknowledgement}

Part of this work was made during a postdoctoral position of the second author in the Department of Mathematics at University of British Columbia (under the supervision of Mike Bennett) and also during
a very enjoyable visit of him to IMB Universit\' e de Bordeaux. He also acknowledges the support of the French-Brazil network.

\end{document}